\documentclass[12pt,draft]{amsart}
\usepackage{amssymb}

\begin{document}

\newtheoremstyle{hplain}%
  {\topsep}
  {\topsep}
  {\itshape}
  {}
  {\bfseries}
  {.}
  { }
  {\thmname{#2 }\thmnumber{#1}\thmnote{ \rm(#3)}}%
  
\newtheoremstyle{hdefinition}
  {\topsep}%
  {\topsep}%
  {\upshape}
  {}%
  {\bfseries}%
  {.}
  { }%
  {\thmname{#2 }\thmnumber{#1}\thmnote{ \rm(#3)}}%
 
\theoremstyle{hplain}
\newtheorem{thm}{Theorem}
\newtheorem{note}[thm]{Note}
\newtheorem{lemma}[thm]{Lemma}
\newtheorem{remark}[thm]{Remark}
\newtheorem{corollary}[thm]{Corollary}
\newtheorem{claim}{Claim}[thm]
\newtheorem{subclaim}{Subclaim}[claim]
\newtheorem{conjecture}[thm]{Conjecture}

\theoremstyle{hdefinition}
\newtheorem{definition}[thm]{Definition}
\newtheorem{question}[thm]{Question}

\newcommand{\restr}{\restriction}
\newcommand{\arr}{\longrightarrow}
\newcommand{\sub}{\subseteq}
\newcommand{\stick}{{}\sp {\bullet}\hskip -7.2pt \shortmid }
\newcommand{\oom}{\stackrel{1-1}{\longrightarrow}}
\newcommand{\iso}{\stackrel{\sim}{\longrightarrow}}
\newcommand{\conc}{{}^\smallfrown}
\newcommand{\setm}{\setminus}
\newcommand{\wprec}{\prec_{\omega_1}}
\newcommand{\wsucc}{\succ_{\omega_1}}
\newcommand{\Hx}[1]{{{\rm H}_{\aleph_{#1}}}}
\newcommand{\Hk}[1]{{{\rm H}_{#1}}}
\newcommand{\Ht}{{\rm H}_\theta}
\newcommand{\mbb}{\mathbb}
\newcommand{\mcal}{\mathcal}
\newcommand{\supf}{\sup{}\!^{(f)}}
\newcommand{\supb}{\sup{}\!^{(b)}}
\newcommand{\Nm}{{\rm Nm}}
\newcommand{\levt}{\tau}
\newcommand{\tg}{{\rm tag}}
\newcommand{\hgt}{{\mathrm{ht}}}
\renewcommand{\frak}{\mathfrak}
\renewcommand{\theequation}{\thesection.\arabic{equation}}

\numberwithin{equation}{section}

\newcommand{\rem}[1]{}
\newcommand{\comment}[1]{\marginpar{\scriptsize\raggedright #1}}

\title{Kurepa-trees and Namba forcing}

\author{Bernhard K\"onig}
\address{\newline
  Boise State University\newline
  Department of Mathematics\newline
  Boise ID 83725-1555\newline
  USA}
\email{bkoenig@diamond.boisestate.edu}

\author{Yasuo Yoshinobu}
\address{\newline
  Graduate School of Information Science\newline
  Nagoya University\newline
  Furocho, Chikusa-ku, Nagoya 464-8601\newline
  Japan}
\email{yosinobu@math.nagoya-u.ac.jp}


\begin{abstract}
  We show that compact cardinals and {\rm MM} are sensitive to
  $\lambda$-closed forcings for arbitrarily large $\lambda$. This is
  done by adding `regressive' $\lambda$-Kurepa-trees in either case.
  We argue that the destruction of regressive Kurepa-trees with {\rm
    MM} requires the use of Namba forcing.
\end{abstract}

\maketitle

\section{Introduction}

Say that a tree $T$ of height $\lambda$ is {\em $\gamma$-regressive}
if for all limit ordinals $\alpha<\lambda$ with ${\rm
  cf}(\alpha)<\gamma$ there is a function $f_\alpha:T_\alpha\arr
T_{<\alpha}$ which is {\em regressive}, i.e.  $f_\alpha(x)<_Tx$ for
all $x\in T_\alpha$ and if $x,y\in T_\alpha$ are distinct then
$f_\alpha(x)$ or $f_\alpha(y)$ is strictly above the meet of $x$ and
$y$. We give a summary of the main results of this paper:

\setcounter{thm}{4}
\begin{thm}\label{intro2}
  For all uncountable regular $\lambda$ there is a $\lambda$-closed
  forcing $\mcal{K}_{\rm reg}^\lambda$ that adds a
  $\lambda$-regressive $\lambda$-Kurepa-tree.
\end{thm}

This is contrasted in Section \ref{compact-sec}:

\setcounter{thm}{6}
\begin{thm}\label{intro3}
  Assume that $\kappa$ is a compact cardinal and $\lambda\geq\kappa$
  is regular. Then there are no $\kappa$-regressive
  $\lambda$-Kurepa-trees.
\end{thm}

Theorems \ref{intro2} and \ref{intro3} establish that compact
cardinals are sensitive to $\lambda$-closed forcings for arbitrarily
large $\lambda$. This should be compared with the well-known result
that a supercompact cardinal $\kappa$ can be made indestructible by
$\kappa$-directed-closed forcings \cite{laver78:_makin}. These results
drive a major wedge between the notions of $\lambda$-closed and
$\lambda$-directed-closed. Another contrasting known result is that a
strong cardinal $\kappa$ can be made indestructible by
$\kappa^+$-closed forcings \cite{gitik89:_hajnal}.
In Section \ref{largerht} we prove

\setcounter{thm}{12}
\begin{thm}\label{intro1}
  Under {\rm MM}, there are no $\omega_1$-regressive
  $\lambda$-Kurepa-trees for any uncountable regular $\lambda$.
\end{thm}
\setcounter{thm}{0}

This shows that {\rm MM} is sensitive to $\lambda$-closed forcings for
arbitrarily large $\lambda$, thus answering a question from both
\cite{koenig04:_fragm_maxim} and \cite{larson00:_separ}. Note that
{\rm MM} is indestructible by $\omega_2$-directed-closed forcings
\cite{larson00:_separ}, so again we find a remarkable gap between the
notions of $\omega_2$-closed and $\omega_2$-directed-closed.
Interestingly enough though, $\omega_2$-closed forcings can only
violate a very small fragment of {\rm MM}. To see this, let us denote
by $\Gamma_{\rm cov}$ the class of posets that preserve stationary
subsets of $\omega_1$ and have the {\em covering property}, i.e. every
countable set of ordinals in the extension can be covered by a
countable set in the ground model. Then we have the following result
from \cite[p.302]{koenig04:_fragm_maxim}:

\begin{thm}
  The axioms ${\rm PFA}$, ${\rm MA}(\Gamma_{\rm cov})$ and ${\rm
    MA}^+(\Gamma_{\rm cov})$ are all indestructible by
  $\omega_2$-closed forcings respectively.\footnote{See below for a
    definition of the axioms ${\rm MA}(\Gamma)$ and ${\rm
      MA}^+(\Gamma)$.}
\end{thm}

So Theorem \ref{intro2} gives

\begin{corollary}\label{regKcons}
  If $\lambda\geq\omega_2$ is regular, then ${\rm MA}^+(\Gamma_{\rm
    cov})$ is consistent with the existence of a $\lambda$-regressive
  $\lambda$-Kurepa-tree.
\end{corollary}

Again, compare this with Theorem \ref{intro1}. It is interesting to
add that ${\rm MA}^+(\Gamma_{\rm cov})$ in particular implies the
axioms ${\rm PFA}^+$ and ${\rm SRP}$. The typical example of a forcing
that preserves stationary subsets of $\omega_1$ but does not have the
covering property is Namba forcing and the proofs confirm that Namba
forcing plays a crucial role in this context. It has already been
established in \cite{larson00:_t} and \cite{MR1713438} that ${\rm
  MA}(\Gamma_{\rm cov})$ can be preserved in an
$(\omega_1,\infty)$-distributive forcing extension in which the
Namba-fragment of {\rm MM} fails. In our case though, the failure of
{\rm MM} is obtained with a considerably milder forcing, i.e.
$\lambda$-closed for arbitrarily large $\lambda$.

The authors would like to thank Yoshihiro Abe, Tadatoshi Miyamoto and
Justin Moore for their helpful comments.

\smallskip The reader requires a strong background in set-theoretic
forcing, a good prerequisite would be \cite{jech97:_set_theor}. We
give some definitions that might not be in this last reference or
because we defined them in a slightly different fashion. If $\Gamma$
is a class of posets then ${\rm MA}(\Gamma)$ denotes the statement
that whenever $\mcal{P}\in\Gamma$ and $D_\xi\;(\xi<\omega_1)$ is a
collection of dense subsets of $\mcal{P}$ then there exists a filter
$G$ on $\mcal{P}$ such that $D_\xi\cap G\not=\emptyset$ for all
$\xi<\omega_1$. The stronger ${\rm MA}^+(\Gamma)$ denotes the
statement that whenever $\mcal{P}\in\Gamma$, $D_\xi\;(\xi<\omega_1)$
are dense subsets of $\mcal{P}$, and $\dot{S}$ is a $\mcal{P}$-name
such that $$\Vdash_{\mcal{P}}\mbox{$\dot{S}$ is stationary in
  $\omega_1$}$$
then there exists a filter $G$ on $\mcal{P}$ such that
$D_\xi\cap G\not=\emptyset$ for all $\xi<\omega_1$, and
$$\dot{S}[G]=\{\gamma<\omega_1:\exists q\in G(q\Vdash_{\mcal{P}}
\check{\gamma}\in\dot{S})\}$$
is stationary in $\omega_1$. In
particular, {\rm PFA} is ${\rm MA}(\mbox{proper})$ and ${\rm MM}$ is
${\rm MA}(\mbox{preserving stationary subsets of $\omega_1$})$. The
interested reader is referred to
\cite{baumgartner84:_applic_proper_forcin_axiom} and \cite{foreman88}
for the history of these {\em forcing axioms}.
  
A partial order is {\em $\lambda$-closed} if it is closed under
descending chains of length less than $\lambda$. It is {\em
  $\lambda$-directed-closed} if it is closed under directed subsets of
size less than $\lambda$. \cite{koenig04:_fragm_maxim} proves that
{\rm PFA} is preserved by $\omega_2$-closed forcings and
\cite{larson00:_separ} that {\rm MM} is preserved by
$\omega_2$-directed-closed forcings.
  
{\em Namba forcing} is denoted by $\Nm$: conditions are trees
$t\sub\omega_2^{<\omega}$ with a {\em trunk ${\rm tr}(t)$} such that
$t$ is linear below ${\rm tr}(t)$ and has splitting $\aleph_2$
everywhere above the trunk. Smaller trees contain more information. It
is known that Namba forcing preserves stationary subsets of
$\omega_1$. If $t\in\Nm$ and $x\in t$ then the last element of $x$ is
also called the {\em tag} of $x$, denoted as $\tg(x)$, and we define
${\rm Suc}_t(x)$ to be the set of tags of all immediate successors of
$x$ in $t$. So ${\rm Suc}_t(x)$ is an unbounded subset of $\omega_2$.
In an abuse of notation, a sequence is sometimes confused with its
tag. We write $[t]$ for the set of infinite branches through $t$.

\section{Stationary limits}

For a tree $T$ and an ordinal $\alpha$, let $T_\alpha$ denote the {\em
  $\alpha$th level of $T$} and
$T_{<\alpha}=\bigcup_{\xi<\alpha}T_\xi$. If $X$ is a set of ordinals,
we write $T\restr X$ for the subtree $\bigcup_{\xi\in X}T_\xi$. The
expression $\hgt(T)$ denotes the {\em height of $T$}. We only consider
trees of functions. If $T$ is a tree and $\mcal{B}$ a collection of
cofinal branches through $T$ then we call $\mcal{B}$ {\em
  non-stationary over} $T$ if there is a function $f:\mcal{B}\arr T$
which is {\em regressive}, i.e.  $f(b)\in b$ for all $b\in\mcal{B}$
and if $b,b'\in\mcal{B}$ are distinct then $f(b)$ or $f(b')$ is
strictly above $b\cap b'$. Otherwise we call $\mcal{B}$ {\em
  stationary over} $T$. A tree $T$ of height $\kappa$ is called {\em
  $\gamma$-regressive} if $T_\alpha$ is non-stationary over
$T_{<\alpha}$ for every limit ordinal $\alpha<\kappa$ of cofinality
less than $\gamma$. The following is easy to check:

\begin{remark}\label{cofinalremark}
  Assume that $A\sub\alpha$ is cofinal in $\alpha$. Then $T_\alpha$ is
  stationary over $T_{<\alpha}$ iff $T_\alpha$ is stationary over
  $T\restr A$.
\end{remark}

The $\omega$-cofinal limits will figure prominently when dealing with
$\omega_1$-regressive trees, so we prove a useful Lemma about these.
For simplicity we only consider trees of height $\omega$. The reader
will notice that the following observations are applicable in Section
\ref{destrKtree}. If $T$ is of height $\omega$ and $\mcal{B}$ a
collection of infinite branches then for any subset $S\sub T$ we let
$$\overline{S}=\{b\in\mcal{B}:b\cap S\mbox{ is infinite}\}.$$
If $S$
is countable and $\overline{S}$ uncountable then we call $S$ a {\em
  Cantor-subtree} of $T$. The class
$\mcal{N}(T,\mcal{B})\sub[H_\theta]^{\aleph_0}$ (for some large enough
regular $\theta$) is defined by letting $N\in\mcal{N}(T,\mcal{B})$ if
and only if there is $b\in\mcal{B}$ such that $b\sub N$ but $b\notin
N$. We have the following

\begin{lemma}\label{statbranchlem}
  Assume that $T$ has height $\omega$ and size $\aleph_1$ and that
  $\mcal{B}$ is a collection of infinite branches. Then the following
  are equivalent:
  \begin{enumerate}
  \item $\mcal{B}$ is stationary over $T$.
  \item
    \begin{enumerate}
    \item Either there is a Cantor-subtree $S\sub T$ or
    \item if we identify $T$ with $\omega_1$ by any enumeration then
      $$E_\mcal{B}=\{\alpha<\omega_1:\sup(b)=\alpha\mbox{ for some
      }b\in\mcal{B}\}$$
      is stationary in $\omega_1$.
    \end{enumerate}
  \item $\mcal{N}(T,\mcal{B})$ is stationary in
    $[H_\theta]^{\aleph_0}$.
  \end{enumerate}
\end{lemma}
\begin{proof}
  The equivalence of (1) and (3) can be found in
  \cite[p.112]{koenig03:_local_coher} and the implication
  $(2)\Longrightarrow(1)$ is easy.
  
  For $(3)\Longrightarrow(2)$, assume $\neg(2)$ and show $\neg(3)$:
  pick an enumeration $e:\omega_1\to T$ such that $E_\mcal{B}$ is
  nonstationary if we identify nodes with countable ordinals via the
  enumeration $e$. Pick a structure $N\prec\Ht$ such that
  $e,T,\mcal{B}\in N$ and set $\gamma=N\cap\omega_1$, so we have
  $\gamma\notin E_\mcal{B}$. Let $b\in\mcal{B}$ be such that $b\sub
  N$. Then $\sup(b)<\gamma$ holds. Now define
  $$\mcal{A}=\{c\in\mcal{B}:\sup(c)=\sup(b)\}.$$
  Note that
  $\mcal{A}\in N$ and $\mcal{A}$ is countable since we know by
  $\neg(2)(a)$ that $\overline{\sup(b)}$ is countable. So
  $\mcal{A}\sub N$, therefore $b\in N$. This shows that
  $N\notin\mcal{N}(T,\mcal{B})$ and $\mcal{N}(T,\mcal{B})$ is
  non-stationary.
\end{proof}

Note that the equivalence of (1) and (3) is to some extent already in
\cite[p.955]{baumgartner84:_applic_proper_forcin_axiom} but our result
differs slightly from this last reference as we have a stronger notion
of non-stationarity. See also \cite{koenig03:_local_coher} for
variations of Lemma \ref{statbranchlem} in uncountable heights.

\section{Creating regressive Kurepa-trees}\label{createKtree}

Let $\lambda$ be a regular uncountable cardinal throughout this
section. We describe the natural forcing $\mcal{K}_{\rm reg}^\lambda$
to add a $\lambda$-regressive $\lambda$-Kurepa-tree and show that this
forcing is $\lambda$-closed.  We may assume the cardinal arithmetic
$2^{<\lambda}=\lambda$, otherwise a preliminary Cohen-subset of
$\lambda$ could be added.  Conditions of $\mcal{K}_{\rm reg}^\lambda$
are pairs $(T,h)$, where
\begin{enumerate}
\item $T$ is a tree of height $\alpha+1$ for some $\alpha<\lambda$ and
  each level has size $<\lambda$.
\item $T$ is $\lambda$-regressive, i.e. if $\xi\leq\alpha$ then
  $T_\xi$ is non-stationary over $T_{<\xi}$.
\item $h:T_\alpha\arr\lambda^+$ is 1-1.
\end{enumerate}

The condition $(T,h)$ is stronger than $(S,g)$ if
\begin{itemize}
\item $S=T\restriction\hgt(S)$.
\item ${\rm rng}(g)\sub{\rm rng}(h)$.
\item $g^{-1}(\nu)\leq_Th^{-1}(\nu)$ for all $\nu\in{\rm rng}(g)$.
\end{itemize}

A generic filter $G$ for $\mcal{K}_{\rm reg}^\lambda$ will produce a
$\lambda$-regressive $\lambda$-tree $T_G$ in the first coordinate and
the sets $$b_\nu=\{x\in T_G:\mbox{there is }(T,h)\in G\mbox{ such that
}h(x)=\nu\}$$
for $\nu<\lambda^+$ form a collection of
$\lambda^+$-many mutually different $\lambda$-branches through the
tree $T_G$. Notice also that the standard arguments for
$\lambda^+$-{\it cc} go through here as we assumed
$2^{<\lambda}=\lambda$.

So we are done once we show that $\mcal{K}_{\rm reg}^\lambda$ is
$\lambda$-closed. To this end, let $(T^\xi,h^\xi)\;(\xi<\gamma)$ be a
descending chain of conditions of length less than $\lambda$. We can
obviously assume that $\gamma$ is a limit ordinal. If the height of
$T^\xi$ is $\alpha^\xi+1$, let
$\alpha^\gamma=\sup_{\xi<\gamma}\alpha^\xi$. We want to extend the
tree $$T^*=\bigcup_{\xi<\gamma}T^\xi,$$
so we have to define the
$\alpha^\gamma$th level: whenever $\nu\in{\rm rng}(h^\xi)$ for some
$\xi<\gamma$, then there is exactly one $\alpha^\gamma$-branch $c_\nu$
that has color $\nu$ on a final segment. Now define
$$T^\gamma_{\alpha^\gamma}=\{c_\nu:\nu\in{\rm rng}(h^\xi)\mbox{ for
  some }\xi<\gamma\}$$
and let $T^\gamma$ be the tree $T^*$ with the
level $T^\gamma_{\alpha^\gamma}$ on top. The 1-1 function
$h^\gamma:T^\gamma_{\alpha^\gamma}\arr\lambda^+$ is defined by letting
$$h^\gamma(c_\nu)=\nu.$$
We claim that $(T^\gamma,h^\gamma)$ is a
condition: the only thing left to check is that
$T^\gamma_{\alpha^\gamma}$ is non-stationary over $T^*$. But this is
witnessed by the function $$f(c_\nu)=\mbox{the $<_T$-least }x\in
c_\nu\mbox{ such that there is }\xi<\gamma\mbox{ with
}h^{\xi}(x)=\nu.$$
Notice that $f$ is regressive: if
$$f(c_\nu)\leq_Tf(c_\mu)\leq_T c_\nu\cap c_\mu,$$
let $\xi$ witness
that $f(c_\mu)=x$, i.e. $h^\xi(x)=\mu$. Then $h^\xi(x)$ must be color
$\nu$ as well since $f(c_\nu)\leq_Tx$ has color $\nu$.  Thus,
$\nu=h^\xi(x)=\mu$.

But $(T^\gamma,h^\gamma)$ extends the chain
$(T^\xi,h^\xi)\;(\xi<\gamma)$, so we just showed
\begin{thm}\label{addlKtree}
  $\mcal{K}_{\rm reg}^\lambda$ is a $\lambda$-closed forcing that adds
  a $\lambda$-regressive $\lambda$-Kurepa-tree.
\end{thm}
We emphasize again that the forcing $\mcal{K}_{\rm reg}^\lambda$ is
{\em not} $\omega_2$-directed-closed but the reader can check that the
usual forcing to add a plain $\lambda$-Kurepa-tree (see e.g.
\cite{jech97:_set_theor}) actually {\em is} $\lambda$-directed-closed.

\section{Destroying regressive Kurepa-trees above a compact cardinal}\label{compact-sec}

If $\lambda$ is a regular uncountable cardinal then a tree $T$ is
called a {\em weak $\lambda$-Kurepa-tree} if
\begin{itemize}
\item $T$ has height $\lambda$,
\item each level has size $\leq\lambda$ and
\item $T$ has $\lambda^+$-many cofinal branches.
\end{itemize}

\begin{lemma}\label{emb-noKtree}
  Suppose that $\lambda$ is a regular uncountable cardinal and there
  is an elementary embedding $j:V\arr M$ such that
  $\eta=\sup(j''\lambda)<j(\lambda)$ and ${\rm cf}^M(\eta)<j(\kappa)$.
  Then there are no $\kappa$-regressive weak $\lambda$-Kurepa-trees.
\end{lemma}
\begin{proof}
  Suppose that $T$ is a $\kappa$-regressive weak $\lambda$-Kurepa-tree
  and $j$ as above. Then there is a regressive function $f_\eta$
  defined on the level $(jT)_\eta$. If $b$ is a cofinal branch through
  $T$, then we find $\alpha_b<\lambda$ such that
  $$f_\eta(jb\restr\eta)\leq_{jT}jb\restr
  j(\alpha_b)=j(b\restr\alpha_b).$$
  Note that if $b$ and $b'$ are two
  distinct branches through $T$ then $jb$ and $jb'$ must disagree
  below $\eta$. Moreover, $j(b\restr\alpha_b)\neq
  j(b'\restr\alpha_{b'})$ holds because $f_\eta$ is regressive. Then
  the assignment $b\longmapsto b\restr\alpha_b$ must be 1-1, which is
  a contradiction to the fact that $T$ has $\lambda^+$-many branches.
\end{proof}

Recall that a cardinal $\kappa$ is {\em $\lambda$-compact} if there is
a fine ultrafilter on $\mcal{P}_\kappa\lambda$. If $\lambda$ is
regular, the elementary embedding $j:V\arr M$ with respect to such a
fine ultrafilter has the following properties:
\begin{itemize}
\item the critical point of $j$ is $\kappa$,
\item there is a discontinuity at $\lambda$,
  i.e. $\eta=\sup(j''\lambda)<j(\lambda)$ and
\item ${\rm cf}^M(\eta)<j(\kappa)$.
\end{itemize}
(see \cite[\S22]{higherinf} for more details). A cardinal $\kappa$ is
said to be {\em compact} if it is $\lambda$-compact for all $\lambda$,
so it follows from Lemma \ref{emb-noKtree} and the above definition:

\begin{thm}\label{noKtree-compact}
  Assume that $\kappa$ is a compact cardinal and $\lambda\geq\kappa$
  is regular. Then there are no $\kappa$-regressive weak
  $\lambda$-Kurepa-trees.
\end{thm}

Using Theorem \ref{addlKtree}, we have

\begin{corollary}
  Compact cardinals are sensitive to $\lambda$-closed forcings for
  arbitrarily large $\lambda$.
\end{corollary}

It was known before that adding a slim\footnote{A $\kappa$-Kurepa-tree
  $T$ is called {\em slim} if $|T_\alpha|\leq|\alpha|$ for all
  $\alpha<\kappa$.} $\kappa$-Kurepa-tree destroys the ineffability of
$\kappa$ and that slim $\kappa$-Kurepa-trees can be added with
$\kappa$-closed forcing.  But note that our notion of regressive is
more universal: slim Kurepa-trees can exist above compact or even
supercompact cardinals.

\section{Oscillating branches}

Now assume that $T$ is an $\omega_2$-tree: we enumerate each level by
letting
\begin{equation}
  \label{eq:level-enum}
  T_\alpha=\{\levt(\alpha,\xi):\xi<\omega_1\}\mbox{ for all }
  \alpha<\omega_2.
\end{equation}
In this situation we identify branches with functions from $\omega_2$
to $\omega_1$ that are induced by the enumerations of the levels. If
$A\sub\omega_2$ is unbounded and $b:\omega_2\arr\omega_1$ is an
$\omega_2$-branch through $T$ then we say that $b$ {\em oscillates on}
$A$ if for all $\alpha<\omega_2$ and all $\zeta<\omega_1$ there is
$\beta>\alpha$ in $A$ and $\xi>\zeta$ such that $b(\beta)=\xi$.

\begin{lemma}\label{osclem}
  Assume that $T$ is an $\omega_2$-Kurepa-tree with an enumeration
  $\levt(\alpha,\xi)\;(\alpha<\omega_2,\,\xi<\omega_1)$ as in
  (\ref{eq:level-enum}) and $A_\iota\;(\iota<\omega_2)$ are
  $\aleph_2$-many unbounded subsets of $\omega_2$. Then there is an
  $\omega_2$-branch $b$ through $T$ that oscillates on every
  $A_\iota\;(\iota<\omega_2)$.
\end{lemma}
\begin{proof}
  Assume not, then for every $\omega_2$-branch $b$ there is
  $\iota_b<\omega_2$ and there are $\alpha_b<\omega_2$,
  $\zeta_b<\omega_1$ such that
  $$b\restr(A_{\iota_b}\setm\alpha_b)\sub\{\levt(\alpha,\xi):
  \alpha\in A_{\iota_b}\setm\alpha_b,\,\xi<\zeta_b\}.$$
  By a cardinality argument we can
  find $\aleph_3$-many branches $b$ such that $\iota_0=\iota_b$,
  $\alpha_0=\alpha_b$ and $\zeta_0=\zeta_b$. But then each of these
  branches is a different branch through the tree
  $$T_0=\{\levt(\alpha,\xi):\alpha\in
  A_{\iota_0}\setm\alpha_0,\,\xi<\zeta_0\}.$$
  $T_0$ has countable
  levels but $\aleph_3$-many branches, a contradiction.
\end{proof}

\section{Destroying regressive Kurepa-trees with {\rm
    MM}}\label{destrKtree}

We introduce a simplified notation for the following arguments: if
$f:t\arr\omega_1$ for some $t\in\mathrm{Nm}$ and $\pi\in[t]$ then we
let
$$\supf(\pi)=\sup_{n<\omega}f(\pi\restriction n).$$
If $b:\omega_2\arr\omega_1$ is an $\omega_2$-branch and
$x\in\omega_2^{<\omega}$ then $b(x)$ really denotes the countable
ordinal $b(\tg(x))$.

\begin{lemma}\label{statNmlem}
  Assume that $T$ is an $\omega_2$-Kurepa-tree and $\mcal{B}$ is the
  set of branches. Let
  $\levt(\alpha,\xi)\;(\alpha<\omega_2,\,\xi<\omega_1)$ be an enumeration
  as in (\ref{eq:level-enum}). Then in the Namba extension $V^{\rm
    Nm}$ there is a sequence
  $$\Delta_G=\langle\delta^G_n:n<\omega\rangle$$
  cofinal in
  $\omega_2^V$ such that
  $$\dot{E}_\mcal{B}=\{\supb(\Delta_G):b\in\mcal{B}\}$$
  is stationary
  relative to every stationary $S\sub\omega_1$ in $V$, i.e.
  $\dot{E}_\mcal{B}\cap S$ is stationary for all stationary
  $S\sub\omega_1$ in the ground model.
\end{lemma}
\begin{proof}
  Assume that $\dot{C}$ is an ${\rm Nm}$-name for a club in
  $\omega_1$, $S\sub\omega_1$ is a stationary set in $V$ and $t_0$ a
  condition in ${\rm Nm}$. Our goal is to find a condition $t_3\leq
  t_0$ and an ordinal $\xi_0\in S$ such that
  $t_3\Vdash\xi_0\in\dot{C}\cap\dot{E_\mcal{B}}$. By a fusion argument
  similar to the ones in \cite[p.188]{MR1713438}, we construct a
  condition $t_1\leq t_0$ and a coloring $f:t_1\arr\omega_1$ such that
  \begin{enumerate}
  \item $f$ increases on chains, i.e. if $v\subsetneq x$ are elements
  of $t_1$ then $f(v)<f(x)$.\label{f-inc}
  \item if the height of $x$ in $t_1$ is odd and $x$ is above the
    trunk then there is $\zeta<\omega_1$ such that
    $$|\{x\conc\beta\in t_1\mid
    f(x\conc\beta)=\xi\}|=\aleph_2\mbox{ for all }\xi>\zeta,$$
    i.e. each ordinal in a final segment of $\omega_1$ has
    $\aleph_2$-many preimages in the set ${\rm
      Suc}_{t_1}(x)$.\label{preimg}
  \item if $G\sub\Nm$ is generic with $t_1\in G$ and
    $\pi:\omega\arr\omega_2^V$ is the corresponding Namba-sequence
    then $\supf(\pi)\in\dot{C}[G]$.\label{sup-C}
   \end{enumerate}
   Given the condition $t_1$, we apply Lemma \ref{osclem} to find a
   branch $b$ that oscillates on all sets ${\rm Suc}_{t_1}(x)\;(x\in
   t_1)$.  Using (\ref{f-inc}) and (\ref{preimg}), we thin out again
   to get a condition $t_2\leq t_1$ with the following property:
   \begin{enumerate}\setcounter{enumi}{3}
   \item if $v\subsetneq x\subsetneq y$ is a chain in $t_2$ above the
   trunk and the height of $x$ is odd then
   $f(v)<b(x)<f(y)$.\label{intertwined}
   \end{enumerate}
   Note that in particular (\ref{preimg}) can be preserved by passing
   to the condition $t_2$, so we may assume that $t_2$ has properties
   (\ref{f-inc})-(\ref{intertwined}). Let us also assume for
   notational simplicity that the height of ${\rm tr}(t_2)$ is even.
   The next step is to find $t_3\leq t_2$ and $\xi_0\in S$ such that
   \begin{enumerate}\setcounter{enumi}{4}
   \item $\supf(\pi)=\xi_0$ for all branches $\pi$ in
     $[t_3]$.\label{inS-cond}
   \end{enumerate}
   To find $t_3$ and $\xi_0$, we define a game $\mbb{G}(\gamma)$ for
   every limit $\gamma<\omega_1$. Fix a ladder sequence
   $l(\gamma)=(\gamma_n:n<\omega)$ for each such $\gamma$. The game
   $\mbb{G}(\gamma)$ is played as follows:
   
   \medskip
  \begin{tabular}{c|ccccc}
       I & $\alpha_0$ & $\alpha_1$ & $\alpha_2$ & $\alpha_3$ & 
       $\ldots$\\\hline
      II & $\qquad\beta_0$ & $\qquad\beta_1$ & $\qquad\beta_2$ & 
      $\qquad\beta_3$ & $\qquad\ldots$
  \end{tabular}
  \medskip\newline where for all $n<\omega$
  \begin{itemize}
  \item $\alpha_n<\beta_n<\omega_2$,
  \item $s_n={\rm tr}(t_2)\conc(\beta_i:i\leq n)\in t_2$ and
  \item $f(s_n)\in(\gamma_n,\gamma)$ whenever $n$ is even.
  \end{itemize}
  II wins if he can make legal moves at each step, so the game is
  determined.
  \begin{claim}
    II wins $\mbb{G}(\gamma)$ for club many $\gamma's$.
  \end{claim}
  \begin{proof}
    Assume not, then there is a stationary $U\sub\omega_1$ such that
    player I wins $\mbb{G}(\gamma)$ for each $\gamma\in U$ via the
    strategy $\sigma_\gamma$. Now pick a countable elementary $N$ such
    that $\xi=N\cap\omega_1\in U$ and $t_2,f,l,U\in N$.
    
    A ladder sequence $l(\xi)=(\xi_n:n<\omega)$ converging to $\xi$ is
    given and we define a sequence $(\beta_n:n<\omega)$ inductively as
    follows: let $\beta_n$ be the least $$\beta>\sup_{\gamma\in
      U}\sigma_\gamma(\beta_i:i<n)$$
    such that
    \begin{itemize}
    \item $s={\rm tr}(t_2)\conc(\beta_i:i<n)\conc\beta\in t_2$ and
    \item $f(s)\in(\xi_n,\xi)$ whenever $n$ is even.
    \end{itemize}
    Such a $\beta$ exists in $N$ by (\ref{preimg}) and elementarity.
    Note that $(\beta_n:n<\omega)$ is a possible record of moves for
    player II if player I goes along with the strategy $\sigma_\xi$.
    But II obviously wins the game $\mbb{G}(\xi)$ if the sequence
    $(\beta_n:n<\omega)$ is played, a contradiction. This proves the
    claim.
  \end{proof}

  Given the claim, pick $\xi_0\in S$ above all $b(x)\;(x\sub{\rm
  tr}(t_2))$ such that II wins the game $\mbb{G}(\xi_0)$. Now we can
  easily find a condition $t_3\leq t_2$ with property
  (\ref{inS-cond}).
  
  If we fix a generic $G\sub\Nm$ with $t_3\in G$ and let
  $\pi_G:\omega\arr\omega_2^V$ be the corresponding Namba-sequence, we
  can define $\delta^G_n=\pi_G(2n+1)$ and
  $\Delta_G=\langle\delta^G_n:n<\omega\rangle$. Then we have
  \begin{enumerate}\setcounter{enumi}{5}
  \item $\supf(\Delta_G)\in\dot{C}[G]$ by (\ref{sup-C}),
  \item $\supf(\Delta_G)=\xi_0$ by (\ref{inS-cond}) and
  \item $\supf(\Delta_G)=\sup_{n<\omega}b(\delta^G_n)=\supb(\Delta_G)$
    by (\ref{intertwined}).
  \end{enumerate}
  But this finishes the proof since
  $$\xi_0\in\dot{C}[G]\cap\dot{E}_\mcal{B}[G]\cap S.$$
\end{proof}

\begin{corollary}
  Assume that $T$ is an $\omega_2$-Kurepa-tree and $\mcal{B}$ the set
  of ground model branches through $T$. Then $\mcal{B}$ is stationary
  over $T$ in the Namba extension.
\end{corollary}

Finally we get the main result for $\omega_2$. We will prove a more
general version of this in Theorem \ref{noreglK}.

\begin{thm}\label{noregK}
  There are no $\omega_1$-regressive $\omega_2$-Kurepa-trees under
  {\rm MM}.
\end{thm}
\begin{proof}
  Assume that $T$ is an $\omega_1$-regressive $\omega_2$-Kurepa-tree
  and that $$\levt(\alpha,\xi)\;(\alpha<\omega_2,\,\xi<\omega_1)$$
  is
  an enumeration as in (\ref{eq:level-enum}). Look at the iteration
  $\mbb{P}=\Nm*{\rm CS}(\dot{E}_\mcal{B})$, where ${\rm
    CS}(\dot{E}_\mcal{B})$ shoots a club through the set
  $\dot{E}_\mcal{B}$ from the statement of Lemma \ref{statNmlem}. The
  poset $\mbb{P}$ preserves stationary subsets of $\omega_1$ by the
  fact that $\dot{E}_\mcal{B}$ is stationary relative to every
  stationary set in $V$. But we have that $\dot{E}_\mcal{B}$ is club
  in $V^\mbb{P}$, so we can use {\rm MM} to get a sequence
  $\Delta=\langle\delta_n:n<\omega\rangle$ converging to
  $\delta<\omega_2$ such that $$\{\supb(\Delta):b\mbox{ is a
    $\delta$-sequence in }T_\delta\}$$
  is club in $\omega_1$. Using
  Lemma \ref{statbranchlem}, we see that $T_\delta$ is definitely
  stationary over $T\restr\Delta$. So $T_\delta$ is stationary over
  $T_{<\delta}$ by Remark \ref{cofinalremark}. Since ${\rm
    cf}(\delta)=\omega$, this contradicts the fact that $T$ is
  $\omega_1$-regressive.
\end{proof}

\section{Larger heights}\label{largerht}

Starting from Theorem \ref{noregK}, we generalize the result to weak
Kurepa-trees in all uncountable regular heights.

\begin{thm}\label{noreglK}
  Under {\rm MM}, there are no $\omega_1$-regressive weak
  $\lambda$-Kurepa-trees for any uncountable regular $\lambda$.
\end{thm}
\begin{proof}
  Since {\rm PFA} destroys weak $\omega_1$-Kurepa-trees (see
  \cite{baumgartner84:_applic_proper_forcin_axiom}), we may assume
  that $\lambda$ is at least $\omega_2$. Now assume that $T$ is an
  $\omega_1$-regressive weak $\lambda$-Kurepa-tree and let
  $\mcal{P}={\rm Col}(\omega_2,\lambda)$ be the usual
  $\omega_2$-directed collapse.  Note that $\mcal{P}$ has the
  $\lambda^+$-{\it cc}, because $\lambda^{\omega_1}=\lambda$ holds
  under {\rm MM} (see \cite{foreman88}). So the tree $T$ has a cofinal
  subtree $T^*$ in $V^\mcal{P}$ that is an $\omega_1$-regressive weak
  $\omega_2$-Kurepa-tree. By throwing away some nodes if necessary, we
  may assume that $T^*$ has the property that
  \begin{equation}
    \label{eq:hom3}
    T^*_x=\{y\in T^*:x\leq_Ty\}\mbox{ has $\aleph_3$-many branches for all
    }x\in T^*.
  \end{equation}
  Now we define an $\omega_2$-directed forcing $\mcal{Q}$ in
  $V^\mcal{P}$ that shoots an actual $\omega_2$-Kurepa-subtree through
  the tree $T^*$: conditions of $\mcal{Q}$ are pairs of the form
  $(S,B)$, where
  \begin{enumerate}
  \item $S$ is a downward-closed subtree of $T^*$ of height $\alpha+1$
    for some ordinal $\alpha<\omega_2$,
  \item $|S|\leq\omega_1$,
  \item $B$ is a nonempty set of branches cofinal in $T^*$ and
    $|B|\leq\aleph_1$,
  \item $b\restr(\alpha+1)\sub S$ for all $b\in B$.
  \end{enumerate}
  We let $(S_0,B_0)\geq_\mcal{Q}(S_1,B_1)$ if $S_0=S_1\restr{\rm
    ht}(S_0)$ and $B_0\sub B_1$.
  
  If $X\sub\mcal{Q}$ is a set of mutually compatible conditions of
  size $\leq\aleph_1$ then we let $S_X$ and $B_X$ be the unions over
  the first respectively second coordinates of $X$. Now $S_X$ can be
  end-extended to a tree $\bar{S}_X$ of successor height by extending
  at least the cofinal branches in the non-empty set $B_X$. But then
  $(\bar{S}_X,B_X)$ is a condition stronger than every condition in
  $X$, hence $\mcal{Q}$ is $\omega_2$-directed-closed. An easy
  cardinality argument shows that $\mcal{Q}$ has the $\aleph_3$-{\it
    cc} because $2^{\aleph_1}=\aleph_2$ holds in $V^\mcal{P}$. It is
  now straightforward that a generic filter $H\sub\mcal{Q}$ will
  produce an $\omega_2$-tree in the first coordinate which is
  $\omega_1$-regressive since it is a subtree of the original tree $T$
  and notice that $\mcal{P}$ and $\mcal{Q}$ both preserve uncountable
  cofinalities. On the other hand, a density argument using
  (\ref{eq:hom3}) shows that the set $$\frak{B}=\bigcup\{B:\mbox{there
    is $S$ such that }(S,B)\in H\}$$
  has cardinality $\aleph_3$, so
  $H$ induces an $\omega_1$-regressive $\omega_2$-Kurepa-tree. The
  composition of two $\omega_2$-directed-closed forcings is again
  $\omega_2$-directed-closed and it was mentioned in the introduction
  that $\omega_2$-directed-closed for\-cings preserve {\rm MM}, so we
  have the situation:
  \begin{itemize}
  \item $V^{\mcal{P}*\mcal{Q}}\models{\rm MM}$
  \item $V^{\mcal{P}*\mcal{Q}}\models$``there is an
    $\omega_1$-regressive $\omega_2$-Kurepa-tree.''
  \end{itemize}
  But this contradicts Theorem \ref{noregK}.
\end{proof}

Using Theorem \ref{addlKtree}, we have

\begin{corollary}
  {\rm MM} is sensitive to $\lambda$-closed forcing algebras for
  arbitrarily large $\lambda$.
\end{corollary}

\bibliography{locco}\bibliographystyle{plain}

\end{document}